\documentclass{article}
\usepackage[utf8]{inputenc}
\usepackage[colorinlistoftodos]{todonotes}

\usepackage{hyperref}
\usepackage{amsmath,amsthm,amsfonts,amssymb,comment,enumerate,wasysym}

\newtheorem*{keyidea}{Key Idea}
\newtheorem*{organizing}{Organizing Principle}
\newtheorem*{metaphor}{Metaphor}

\title{User's Guide Project: Looking Back and Looking Forward}
\author{Don Larson, Kristen Mazur, David White, Carolyn Yarnall}
\date{June 2019}

\begin{document}

\maketitle
\begin{abstract} 

In 2014 Luke Wolcott created the User's Guide Project in which a group of algebraic topologists came together to write user's guides to coincide with their research papers in hopes of making their research more accessible. We examine the role of this innovative project within the greater mathematics community. We discuss the structure and history of the project, its impact on the community, and its value to the participants of the project. We end by encouraging the math community to recognize the value of the project and expand the User's Guide Project to other subfields.
\end{abstract}

\section{Introduction}

Over the past 100 years mathematics has become increasingly specialized, and currently it can be difficult to understand  research in fields other than your own. This specialization also means that non-mathematicians often have no concept of what mathematical research entails. Further, most academic mathematicians are judged by their publication records: both the quantity of papers and the prestige of the journals they appear in. Consequently, mathematicians feel compelled to write in a very specific style, aiming to present the results in the most concise language possible and to impress the reader with their brilliance.

The User's Guide Project grew from Luke Wolcott's desire  to push back against these influences in order to make math research more human and accessible. He envisioned a world in which every research paper came paired with another, less formal user's guide that not only explained the big picture ideas of the paper but also discussed how the author(s) developed their story \cite{wolcott-jhm}. Hence, over a span of four years, twelve mathematicians wrote user's guides to accompany one of their research papers. In these guides, the mathematicians wrote for understanding (rather than concision), summarized their research in non-mathematical terms,  explored the humanistic aspects of the experience of conducting mathematical research, and revealed the true story (with all its trials and tribulations) of how the associated research paper came to be.  In short, the user's guides represented an entirely different world of  writing that focused on mathematical play and personal experience instead of technical results.

We are not the first to recognize the importance of play in mathematics. In \cite{su}, Francis Su states,  ``\textit{mathematical play builds virtues that enable us to flourish in every area of our lives\dots math play builds hopefulness\dots  math play builds community\dots math play builds perseverance.}” However, the play that brings us to our research is often lost in our published works. For example, it was play that drew the second author to the concept of Mackey functors. Yet, in order to publish her paper, \textit{An Equivariant Tensor Product on Mackey Functors} \cite{hill-mazur}, it needed to be as succinct as possible while including all pertinent mathematical detail (and not one thing more). Hence, no semblance of play remained in the final version. Moreover, only readers that are fluent in the language of category theory and have advanced knowledge of algebra and algebraic topology can fully understand her paper.  However, the beauty and play that made her fall in love with Mackey functors were not lost because she wrote a user's guide to accompany the technical paper.  

Furthermore, user's guides  draw attention to important ideas that are  hidden in small corners of unpublished documents  or that are only shared orally at conferences and seminars. The key idea that the first author used for his paper, \textit{The Adams-Novikov $E_2$-term for Behren's Spectrum $Q(2)$ at the Prime 3} \cite{larsonthesis}, came from an informal note listed in an ``Other Documents" section of a leading researcher's personal webpage. This idea is now explicitly stated as a Key Idea in the corresponding user's guide \cite{larson}.  Similarly,  the fourth author's user's guide \cite{yarnall} includes diagrams of her constructions that she shared at conferences but were not suitable for a published article.

After three volumes and 13 guides, the User’s Guide Project is currently dormant as Luke (who was the head editor as well as the founder of the project) has moved on from academia. Now four of the participants have come together to evaluate the impact of the User's Guide Project. In the sections that follow, after providing an overview of the project, we discuss the value of the User's Guide Project to the mathematics community, the value of the project to its participants, and the project’s future in mathematics.

\section{Content,  process, and origin} \label{sec:history-process}

The User's Guide Project ran from 2014 until 2017. Twelve authors wrote a structured user's guide to accompany one of their research papers.  A description of the User's Guide Project was published in 2015 \cite{humanistic} and is summarized in Section \ref{subsec:topics} below, along with numerous examples not present in \cite{humanistic}.  
Full details on the User's Guide Project can be found on the project website \cite{enchiridion}. In Section \ref{subsec:process}, we discuss the process of writing a user's guide, and in Section \ref{subsec:history}, we discuss the origin of the User's Guide Project.

\subsection{User's guide content} \label{subsec:topics}

Each user's guide consists of four sections, or topics that are devoted to ``key insights'' of the paper, the ``metaphors and imagery'' the authors associate to the research, the ``story of the development'' of the research, and a ``colloquial summary'' of the research. We provide further details on each topic below.

\paragraph{Topic 1: Key insights and central organizing principles.} Topic 1  functions like an introduction but with more emphasis on concepts and exposition compared to the introduction of the corresponding research paper. The goal is to help the reader learn, rather than to impress a referee or editor. In this section, key ideas and organizing principles are stated and numbered, much like theorems  in published papers. For example, one user's guide \cite{white} states 

\begin{organizing}
Whenever possible, one should work in a setting where it is possible to replace the objects of interest by nicer objects which are equivalent in a suitable sense.
\end{organizing}

Another user's guide \cite{frankland} states 

\begin{keyidea}
Reduce statements about $L$-complete $E_*$-modules to statements modulo the maximal ideal $\mathfrak{m} \subset E_*$. 
\end{keyidea}

As the reader can see, these ideas and principles can vary widely in terms of how ``mathy'' they are. Topic 1 often also includes the main theorems in the associated research paper and gives a general overview of the structure of the proofs of these theorems.

\paragraph{Topic 2: Metaphors and imagery.} In Topic 2 the authors describe the ``right way to think about'' \cite{humanistic} their research using the types of helpful descriptions, mental imagery, and conceptual metaphors often given in a seminar or conference talk. It is very common for this topic to contain pictures and displayed diagrams that often do not make it into the published version of the research paper because of the community emphasis on concision. For instance, the user's guide \cite{mazur} contains six images spread over five pages. The user's guide \cite{yarnall} contains ten figures and tables spread over seven pages.  Colorful language is also common, e.g. one user's guide \cite{larson} uses the metaphor of a hotel with infinitely many levels to explain the complexity present within the stable homotopy groups of spheres. 
Another example of a metaphor from Topic 2 is given below.  

\begin{metaphor}\cite{lorman}
Slice  up  the  given  object  of  interest  into  more  manageable pieces, compute those pieces, then put them back together.
\end{metaphor}

Topic 2 also allows authors an opportunity for introspection: to tug at the strings of intuition that led to their research results, to think about {\em why} they chose their research, and to think about the imagery associated with the research process itself. For instance, the user's guide \cite{white} presents a metaphor comparing the choice to work with model categories instead of $\infty$-categories to a preference for detail in art over impressionist paintings. 

\paragraph{Topic 3: Story of development.} Topic 3 is the most human aspect of the guide, and in it authors give an honest and personal account of their research process. This can help the reader understand the importance of conversations with an advisor, with office mates, on hikes during conferences, in bars after a seminar, etc. Moreover, many user's guides use this section to relate stories of setbacks. Consequently, this section allows the reader to understand that setbacks are common occurrences in research, that  the first proof is rarely the best proof, and that the development of a paper is often far from linear. 

For example, one author \cite{beardsley} describes ``moderately frantic conversations with my advisor and my de facto co-advisor.'' Another \cite{yeakel} wrote about her dismay at discovering an error in her results after having already presented them at a conference. Another \cite{malkiewich} explained the research process by saying, ``the process was very slow \dots  I would stare at incomprehensible papers, make laughably n\"{a}ive guesses as to what was going on, prove the guesses were wrong, make slightly less wrong guesses, and continue.'' 

Sometimes, authors use Topic 3 to discuss when, where, and how mathematical breakthroughs occurred. For many, breakthroughs happened when visiting researchers at other universities \cite{larson}, during visits to give a seminar talk \cite{malkiewich}, or at a conference \cite{frankland}. For some \cite{wolcott} breakthroughs happened while riding the trans-Siberian railway or during a motorcycle trip from the south of France to the Alps \cite{wolcott2}. Other authors used Topic 3 to write about the lengthy experience of revising a thesis and turning it into a publishable paper \cite{mazur}, or about the publication time-line from submission, to referee report, to resubmission, to eventual acceptance of the paper \cite{frankland}. 

Lastly, some authors used Topic 3 to relate deeply personal stories. For instance, one author \cite{white} wrote about how his ability to conduct research was hampered  by the illness of his father and how he eventually got back to math after learning important coping techniques. We feel that graduate students and early career researchers benefit tremendously from reading about the personal experiences of the user's guide authors, and we discuss these benefits further in Section \ref{sec:value-to-community}.

\paragraph{Topic 4: Colloquial summary.} The final topic of a user's guide gives a colloquial summary of the corresponding research paper, written for an entirely non-mathematical audience, or in some cases written for ``the Curious Undergraduate'' \cite{mazur}. This section is meant to demystify what mathematicians actually do, the ways we think about the objects we study, and our passion for our work. In this section, it is common to find the phrases many mathematicians use to describe their work to family and friends, e.g. comparing donuts and coffee mugs \cite{frankland}, explaining electrical flows on topological spaces by way of analogy with traffic patterns \cite{catanzaro}, or explaining ``in the eyes of a topologist two spaces are deemed equivalent if one can continuously be transformed into the other without breaking it apart'' \cite{merling}. Some authors discuss the experience of conducting mathematical research, e.g. \cite{beardsley} contains the nugget (not often revealed outside of mathematics) that ``mathematicians are blindly groping in the dark, searching for structure, and every once in a while they find the exact same beautiful structure in two different places at once.'' Other user's guides contain metaphors relating mathematics to apple pie \cite{wolcott}, Sesame Street \cite{larson}, and even Taylor Swift lyrics \cite{mazur}.

While the structure of the user's guides was heavily influenced by Luke Wolcott's unique interests and experiences, authors were given leeway in how to interpret the four topic prompts, and this leeway led to the plethora of examples presented above. Giving authors leeway also allowed them to have fun while writing their user's guides, and was an important piece of recruiting authors to write a user's guide at all. For instance, some authors used the experience of writing a user's guide to sharpen their expository writing skills, while others used the experience to explore the humanistic aspects of their own research.

\subsection{Process of writing a user's guide} \label{subsec:process}

Creating a volume of user's guides is a collaborative process. Each cohort (i.e. the authors of a volume) wrote the topics one by one over the course of a year, approximately one topic every four months. Most of the authors found writing a user's guide quite challenging as they had no formal training in the type of writing required for a user's guide. Thus, the writing process involved large amounts of feedback and support from the other authors in the cohort. Each topic was  read by every member of the cohort, who then provided detailed feedback to the author on both its mathematical correctness  and on its exposition. After writing the last topic, the topics were compiled into a full user's guide, and the cohort peer reviewed each others' guides one last time. Throughout the peer review process, the authors focused on contemplation and structured reflection, which, in Luke's words, ``would result in better mathematicians doing better mathematics'' \cite{wolcott-jhm}.

Because many universities only ``count'' research that is peer reviewed, it was important that all user's guides be peer reviewed. However, the type of peer review in the User's Guide Project is somewhat different than traditional peer review. For one thing, the user's guide peer review was not anonymous. It was also much deeper than the peer review one receives from a journal referee, and often encouraged the author to add more examples, metaphors, and imagery, rather than encouraging the author to remove text in favor of concision. In this paper, when we say ``peer review'', we will always mean the deeper, more personal peer review that is a crucial component of the User's Guide Project. 

\subsection{Origin of the project} \label{subsec:history}

The reader may wonder why the user's guides consist of the four topics laid out above. In this section, we answer this question by relating how the User's Guide Project came to be.

The User's Guide Project was the brainchild of Luke Wolcott. Luke received his PhD in homotopy theory in 2012, and worked as a mathematics professor until 2017. His vision for the User's Guide Project stemmed from his interest in humanistic mathematics that dated from at least 2007, when he ran a survey to answer the question ``What does math sound like?" \cite{wolcott-web}. For instance, the topic on key ideas was inspired by Luke's experiences with writing about being a mathematician (\cite{wolcott-blog} and \cite{wolcott-jhm}) and serving as editor and writer for the AMS Graduate Student Blog \cite{ams-blog}. His interest in the metaphors and imagery associated with mathematics was informed by his collaboration with artists to make mathematical art and dance \cite{wolcott-web} and by his published papers on math art \cite{wolcott-art1, wolcott-art2, wolcott-art3} and on math poetry \cite{wolcott-poetry}.

Before beginning the User's Guide Project Luke ran a seminar on the human dimensions of mathematics research \cite{wolcott-web},  blogged about his personal struggle in doing math research \cite{wolcott-blog}, and  published a book connecting his experiences as a mathematician to his  experiences in life outside of mathematics \cite{wolcott-book}. Hence, it was natural for him to develop a user's guide system that included a section in which the authors discussed how they did their research.

Luke's 2012 PhD thesis contained the seeds of the User's Guides Project, via subsections devoted to the ``Experiential Context'' of the results \cite{wolcott-thesis}. These subsections focused on the story of the development of the results and on colloquial summaries. A subsequent paper \cite{wolcott-jhm} by Luke argues for the value of these subsections and can be viewed as a call to create the User's Guide Project. That call was answered in 2014 when Luke teamed up with David White to recruit fellow young algebraic topologists, Cary Malkiewich, Mona Merling, and Carolyn Yarnall to write the first volume of user's guides \cite{enchiridion}.

All of the guides in the three volumes of the User's Guide Project were written by algebraic topologists. However, this does not represent a belief that topology  requires user's guides more than other fields of mathematics and instead was entirely based on the fact that Luke and David had strong networks only in algebraic topology. It is worth noting, however, that algebraic topology is a very visual field, and so algebraic topologists might be especially well-suited to writing about the metaphors and imagery they associate with their research. Regardless, the plan was always to extend the User's Guide Project to other branches of math.  Because of the essential role of peer review in the process of crafting a user's guide, we encourage future user's guide organizers to assemble  cohorts of authors from the same field of mathematics, e.g. a cohort consisting solely of representation theorists rather than  some representation theorists and some complex geometers.

\section{Value to the mathematics community} \label{sec:value-to-community}

The User's Guide Project can help mitigate 
many issues within the math
community, and here we  highlight
three issues that specifically arise in algebraic topology.  
The first  is overly-technical 
mathematical
machinery that makes it hard for fledgling algebraic
topologists to get situated and hinders
those outside the field
from leveraging its tools and accurately
judging topological work. 
The second issue is the predominance of 
algebraic topology ``folklore;'' that is,
knowledge passed along through
conversation that is never written down
carefully.  
The third issue is a culture that encourages
hiding weakness in order to preserve
reputation.  Thankfully, 
writings that endeavor to combat one or more of the aforementioned issues are sprinkled throughout the literature and the mathematical blogosphere (see the list at the
end of this section).  But these writings 
appear all
too infrequently and are ultimately
read only by those willing to comb through the footnotes and appendices of resources filled primarily with mathematical technicalities.  By contrast,
user's guides make such
writings the main attraction.

Consider the following hypothetical situation:
a paper containing a significant result
appears on the arXiv with an abstract that begins 
``Let $E$ be an $E$-infinity ring spectrum.''
This is quite plausible, but some would contend
(e.g., \cite{hovey}) that the first sentence does not bode well for the paper or its readership.  Let us speculate as to why. 
Suppose the paper's author  submitted the work
to a wide-audience
journal, such as Journal of the AMS.  
Suppose also that a graduate student in algebraic
topology and a researcher from
another field 
both see the paper on the arXiv
and suspect the paper might be helpful for their
own work.  Ideally,
the paper would be accepted to the journal, the grad student
would find the paper readable and be able to 
leverage it for her dissertation, and the outside researcher would
quickly surmise the connections with her field
and use them to build a bridge to algebraic topology.  
And yet, it is quite conceivable that none of these things would occur, 
because the concept of an
$E$-infinity ring spectrum is a prime example
of technical mathematical machinery.  
Unless the author is able to quickly bring this 
machinery down to earth 
via intuitions and insights 
that are meaningful to the wide-audience
journal editor, the grad student, and the outside researcher,
the paper is likely to be ignored.

But, what if the author also wrote a user's guide to 
accompany the paper that provided
intuition for how to think about an $E$-infinity ring spectrum in   Topic 2: Metaphors and imagery, like the way 
\cite{larson} does for Greek letter elements in the stable
homotopy groups of spheres?   What if this user's guide
chronicled how these intuitions were obtained in Topic 3: Story of  development, like the way \cite{merling} does for equivariant
bundles?   
Then the guide could come to the rescue 
in a couple of ways.  For one, the grad student and outside researcher could leverage the user's guide to more easily trudge
beyond the author's
daunting opening sentence.  For another,
the author could revamp their abstract using the intuition provided in the user's guide, thereby
making the abstract more inviting and favorable with the wide-audience journal editor.

The heavy reliance on word-of-mouth
for knowledge dissemination in algebraic
topology---a problem the 
User's Guide Project is tailor-made to solve---is
well known.  
For example, in June 2019, a prominent algebraic topologist emailed
a research community 
mailing list because he
could not track down a written reference 
to a spectral sequence that he used ``all the time without thinking.''  This knowledge accessibility problem leads to a lack of inclusiveness.  
What is often required to track down some crucial
algebraic topology fact is a
conversation with the right person, at the right conference, under the right circumstances.  One also has to know the right questions to ask and be sufficiently outgoing to ask them. 
Moreover, one needs to not only
{\em know of} the right person to ask,
but also {\em know} that person sufficiently 
well in order to be comfortable 
making the approach.  
Some among us may never see the planets align like that.
It is well-documented
(see, e.g., \cite{FMP}, \cite{GRD}, \cite{herzig}) that 
the sense of belonging one would need to pursue
information in the manner described above
is often lacking among members of 
under-represented populations, putting them at a
particular
disadvantage.  
User's guides 
put these oral traditions down on paper and make them widely accessible.
By providing a suitable open-access repository for mathematical folklore (that
undergoes 
an intensive and 
collaborative peer-review
process, to boot), the User's Guide
Project is uniquely positioned
to bridge potential knowledge gaps between those with big travel budgets and gregarious personalities, and those without, thereby making algebraic topology and mathematics in general more inclusive.

Young mathematicians often experience large amounts of distress and self-doubt in the process of establishing a research portfolio.  Sir Michael Atiyah once shared that both he and Jean-Pierre Serre felt like quitting research early in their respective careers
\cite{Atiyah}.
Even so, much to the detriment of our mental health, most of us bottle up these feelings
out of a fear that showing such ``weakness" could adversely
affect our careers.  The user's guides, however, provide a forum for sharing distress and self-doubt, assuring others that they are not alone.  Further, many of the user's guides also discuss other feelings, experiences, and teachable moments  such as how to schedule time for research, how to juggle several projects simultaneously, how to find new problems to work on, helpful conversations with colleagues, humorous anecdotes, pursuing dead ends, the elation of making a breakthrough, and the frustrating and/or beneficial interactions with journal referees.  

The issues outlined here, and 
the ways in which the User's Guide Project
addresses them, are by no means exclusive
to the field of algebraic topology.  
Mathematicians across all disciplines
experience the 
sorts of difficulties discussed above. 
Every field has some amount of overly-technical machinery, and every field
has a knowledge
gap.  Lastly and most importantly,
{\em everybody struggles}.  It needs
to be said more.  It needs to be
chronicled more, and the User's Guide Project
provides an ideal place to do it.  

We conclude this section by making a partial list of quality content that 
may be useful in the absence of a current
User's Guide Project.

\begin{enumerate}
\item Omar Camarena's writings (\texttt{https://www.matem.unam.mx/$\sim$omar/})
\item Tim Gowers' blog
(\texttt{https://gowers.wordpress.com/})
\item Living Proof: Stories of Resilience Along the Mathematical Journey \\ (\texttt{https://www.ams.org/about-us/LivingProof.pdf})
\item Fosco Loregian's writings (\texttt{http://www.math.muni.cz/$\sim$loregianf/})
\item Cary Malkiewich's blog (\texttt{https://highlyconnected.home.blog/})
\item Akhil Mathew's writings (\texttt{http://math.uchicago.edu/$\sim$amathew/})
\item Dan Murfet's blog (\texttt{http://therisingsea.org/})
\item Eric Peterson's book {\em Formal Geometry
and Bordism Operations} \cite{peterson}
\item Tim Porter in nLab (\texttt{https://ncatlab.org/nlab/show/Tim+Porter})
\item Neil Strickland's site (\texttt{https://neil-strickland.staff.shef.ac.uk/})
\item Terry Tao's blog
(\texttt{https://terrytao.wordpress.com/})
\item Ravi Vakil's blog (\texttt{math.stanford.edu/$\sim$vakil/216blog/index.html})
\end{enumerate}

In addition to positively impacting the math community, the User's Guide Project greatly benefits the authors of the guides as we now address.

\section{Value to the user's guide authors}\label{sec:authorvalue}

Most authors chose to participate in the User's Guide Project because (1) they agreed with the mission of making algebraic topology more accessible and of humanizing the mathematical research process, and/or (2) they hoped to improve both their expository writing skills and their understanding of their research. Indeed, the benefits of participating in the User's Guide Project include the rare opportunity (especially for recent PhDs) to reflect on the research process and to write a mathematical document that is comprehensible to a wide audience. 

Moreover, another (perhaps unexpected) benefit of participating in this project is the sense of community amongst the cohort of authors of each volume. Despite often being scattered across the globe, the peer-editing component of the project gave authors a sense of belonging, and reading and commenting on the research stories of their peers helped authors embrace their own stories of discovery through struggle. In a survey administered in the summer of 2018\footnote{Survey received exempt status from the Institutional Review Board at Denison University.}, the user's guide authors were asked what they liked  best about  the project. Six of the 11 authors surveyed gave responses centered around the community that developed during the project. For example, one author wrote that they liked, ``getting to share the process of doing mathematics," while another said that the best part of the project was ``realizing that others ran into the same pitfalls as me." Another response stated that their favorite aspect was ``going beyond subjectivity, to co-discover the intersubjective reality of mathematical experience..."

The community created through the User's Guide Project was especially powerful to  the coauthors of this article who wrote their guides while working at teaching-focused institutions. Since none of our colleagues at our institutions were specialists in algebraic topology, we often felt isolated from the algebraic topology community, making it difficult to continue working on our research projects. Thus, the User's Guide Project community became a vital component to our research success by allowing us to stay connected to and hence active in the field.

In addition to developing a strong sense of community, this project improved the participants' writing capabilities since this project created a unique opportunity to write and peer-review non-technical and expository papers. Writing a user's guide required authors to reflect on the big picture of their research and to communicate this picture to a general audience. Writing Topic 4 in particular was excellent practice for writing abstracts for grant proposals since these abstracts often needed to be readable by non-mathematicians.  Further, the peer-review process  enabled authors to refine their writing by learning from giving and receiving feedback. Thus, this exercise not only improved the authors' communication skills, but also improved their understanding of the field of algebraic topology. It is not surprising that when asked about the best part of writing the guide, four of the 11 authors surveyed discussed the fact that it improved their writing skills and allowed them to think about the big picture of their research. One author even wrote, ``I like grappling with the basic ideas underlying my thesis and thinking hard about how to discuss them intuitively, with minimal use of technical machinery."

Furthermore, the User's Guide Project gave the authors an opportunity to reflect on the human aspects of the math research process, from the excitement of discovery to the frustration of failure. While many mathematicians discuss such things on blogs or via social media, the User's Guide Project gave a systematic method for developing a coherent and honest account of the research process. Indeed, the project (especially Topic 3 of the guide) provided a supportive environment that encouraged the authors to write short yet genuine memoirs detailing their personal research experience. (See Section \ref{subsec:topics} for examples.) Many of the authors appreciated this humanistic aspect of the project. In fact, when asked what they liked best about the User’s Guide Project, one author wrote, ``I enjoyed reflecting on how I came to these results and on my growth during graduate school."

In general, the authors of the user's guides agree that participating in the User's Guide Project was a positive experience. Eight out of the 11 authors surveyed  said that the project was valuable for understanding their own research, ten of the authors said that it was a valuable exercise in clear mathematical writing, and seven authors said that writing the user's guide was a valuable use of time. However, many of the authors also worry that their home institutions would not agree that the project is worthwhile. When asked if their current institution values the fact that they wrote the user's guide (in the sense of promotion and tenure) only one author said yes while six authors said no and four were unsure of how the project would be perceived. Moreover, when asked about their least favorite aspect of the guide, most authors discussed the fact that it was a large time commitment that took time away from other “more legitimate” research projects. Yet, despite this fact, eight of the 11 authors surveyed still recommended participating in the User's Guide Project, which indicates that the benefits of this project outweigh the lack of official recognition and institutional support.

\section{Future and call to arms} \label{sec:call-to-arms}

As highlighted in previous sections, the inaccessibility of mathematics is problematic. This inaccessibility makes it difficult for those who are math-curious but not classically math-trained  to learn about the field. Even established mathematicians often have trouble expanding their research into a different subfield of mathematics. The depth of knowledge often required for mathematical exploration can stifle creativity and obscure the beauty and play greatly valued by the coauthors of this article, Francis Su \cite{su}, and the editors and many contributors to this journal (e.g., \cite{huber-karaali} and \cite{storm-zullo}). The User's Guide Project is one way to break the barrier of inaccessibility in order to grow and strengthen the math community.

This inaccessibility is especially troublesome for graduate students, young researchers with temporary employment, and members of underrepresented groups. Providing support for these individuals to engage deeply with mathematical research is essential. The User's Guide Project can help such individuals connect to the community and better understand  previously unobtainable mathematical knowledge. These guides can even provide insight on surviving graduate school (e.g., \cite{mazur} and \cite{white}). The user's guides humanize the authors, comfort readers facing roadblocks in their research, and give readers hope to persevere in their careers as mathematicians. They also contain sage advice, playful anecdotes, and concrete tips for success.   
Moreover, the authors themselves benefited from participating in the User's Guide Project as demonstrated by the survey responses in Section \ref{sec:authorvalue}.

We call upon the math community and institutions of higher learning to value expository
projects like the User's Guide 
Project in the context of hiring, promotion, and tenure. 
It is frustrating that efforts to illuminate mathematics often do not receive official support or recognition, while research articles that only a handful of mathematicians comprehend are valued. While the particular User's Guide Project begun by Luke is currently dormant, the authors feel that it should continue in algebraic topology and spread to other fields. 

Specifically, the authors encourage the reader to begin their own User's Guide Projects. In order to do this, we offer a few helpful resources and suggestions:
\begin{itemize}
\item Collaboration: As discussed in Section \ref{subsec:process}, the importance of having a group of collaborators with which to embark on the user’s guide adventure cannot be understated. The original three cohorts were comprised of (relatively) early career mathematicians. While this is not necessary, several authors were motivated to provide clear exposition as they remembered  struggling to understand difficult concepts when they began their research project. As for the logistics of collaboration, while the first collaborators communicated primarily via email, the second and third groups worked together through CoCalc (previously SageMath). However, with the advent of the online LaTeX editor Overleaf and video conferencing platforms such as WebEx and Zoom, the authors recommend writing the guides via Overleaf and collaborating via a video conferencing platform.

\item Format/Content: We give details on the format and content of the User's Guide Project in Section \ref{subsec:topics}. This format works well though one should feel free to adapt the format to suit their needs. 

\item Sharing/Publication: A simple forum for sharing  user's guides with the math community could be a website as was used for the original project. Further, user's guide cohorts may wish to chronicle the creation of their guides in a blog-style series of posts or publish individual guides in a journal for which more expository work is appropriate. 

\item Potential Obstacles: One of the largest hurdles encountered in beginning such a project is finding time to devote to this creative activity that cannot be classified as traditional research. Some people declined to participate in the project because of this reason. As mentioned at the end of Section \ref{sec:authorvalue}, even some of the user's guide authors found it difficult to balance time spent on this project with time needed for original research publications. Additionally, the guides themselves were difficult to write because the style of writing and topics covered differ from research articles. Finally, the lead task of overseeing the project can be overwhelming. Organizing the collaboration, keeping track of the authors’ progress, and spurring them to action is not for the faint-hearted. Despite these obstacles,  the intrinsic value of participating in such a project makes the difficult writing exercise and reallocation of time worthwhile.
\end{itemize}

We invite any reader who wishes to begin their own User’s Guide Project to contact the coauthors of this paper for further insight and advice.

\end{document}